\documentclass[reqno]{amsart}

\usepackage{amscd,amsmath,amssymb,amsfonts}
\usepackage{stmaryrd}
\usepackage{bbold}
\usepackage{times}
\usepackage[all]{xy}
\usepackage{mathrsfs}
\usepackage{graphicx}

\usepackage{amsthm}

\theoremstyle{plain}
\newtheorem{thm}{Theorem}
\newtheorem{lem}[thm]{Lemma}
\newtheorem{cor}[thm]{Corollary}
\newtheorem{prop}[thm]{Proposition}

\theoremstyle{definition}

\newtheorem{claim}[thm]{Claim}
\newtheorem{example}[thm]{Example}
\newtheorem*{ack}{Acknowledgement}

\newcommand{\Ker}{\operatorname{Ker}}

\newcommand{\Hom}{\operatorname{Hom}}
\newcommand{\End}{\operatorname{End}}
\newcommand{\Aut}{\operatorname{Aut}}

\newcommand{\stab}{\operatorname{stab}}
\newcommand{\ann}{\operatorname{ann}}
\newcommand{\lann}{\operatorname{l.ann}}
\newcommand{\rann}{\operatorname{r.ann}}
\newcommand{\Qr}{\operatorname{Q}_{\text{\rm r}}}
\newcommand{\Ql}{\operatorname{Q}_{\ell}}
\newcommand{\Qcl}{\operatorname{Q}_{\text{\rm cl}}}
\newcommand{\Qmax}{\operatorname{Q}_{\text{\rm max}}}
\newcommand{\cCl}{\mathcal{C}_{\ell}}
\newcommand{\ad}{\operatorname{ad}}

\newcommand{\Spec}{\operatorname{Spec}}
\DeclareMathOperator{\Rat}{Rat}
\DeclareMathOperator{\Prim}{Prim}

\newcommand{\GSpec}{G\text{-}\!\Spec}
\newcommand{\GRat}{G\text{-}\!\Rat}

\newcommand{\Fract}{\operatorname{Fract}}

\newcommand{\GL}{\operatorname{GL}}

\newcommand{\Id}{\operatorname{Id}}

\newcommand{\onto}{\twoheadrightarrow}

\newcommand{\into}{\hookrightarrow}
\newcommand{\tto}{\longrightarrow}
\newcommand{\iso}{\stackrel{\sim}{\tto}}

\renewcommand{\bar}[1]{\overline{#1}}

\newcommand{\til}[1]{\widetilde{#1}}

\newcommand{\Gal}{\operatorname{Gal}}

\renewcommand{\phi}{\varphi}

\renewcommand{\k}{\mathbb{k}}

\newcommand{\cB}{\mathcal{B}}
\newcommand{\sE}{\mathscr{E}}

\newcommand{\cC}{\mathcal{C}}

\newcommand{\Htil}{G_{(P:H)}}
\newcommand{\Itil}{\widetilde{I}}
\newcommand{\fg}{\mathfrak{g}}
\newcommand{\fk}{\mathfrak{k}}

\newcommand{\fm}{\mathfrak{m}}

\newcommand{\Gm}{\mathbb{G}_{\text{\rm m}}}
\newcommand{\Ga}{\mathbb{G}_{\text{\rm a}}}

\newcommand{\ZZ}{\mathbb{Z}}

\newcommand{\cen}{\mathcal{Z}}

\newcommand{\upin}{\text{\ \rotatebox{90}{$\in$}\ }}

\begin{document}




\title[Rational ideals]%
{Group actions and rational ideals}

\author{Martin Lorenz}

\address{Department of Mathematics, Temple University,
    Philadelphia, PA 19122}

\email{lorenz@temple.edu}

\urladdr{http://www.math.temple.edu/$\stackrel{\sim}{\phantom{.}}$lorenz}

\thanks{Research of the author supported in part by NSA Grant H98230-07-1-0008}



\begin{abstract}
We develop the theory of rational ideals for arbitrary associative
algebras $R$ without assuming the standard finiteness conditions,
noetherianness or the Goldie property. The Amitsur-Martindale ring
of quotients replaces the classical ring of quotients which
underlies the previous definition of rational ideals but is not
available in a general setting.

Our main result concerns rational actions of an affine algebraic
group $G$ on $R$. Working over an algebraically closed base field,
we prove an existence and uniqueness result for generic rational
ideals in the sense of Dixmier: for every $G$-rational ideal $I$ of
$R$, the closed subset of the rational spectrum $\Rat R$ that is
defined by $I$ is the closure of a unique $G$-orbit in $\Rat R$.
Under additional Goldie hypotheses, this was established earlier by
M{\oe}glin and Rentschler (in characteristic $0$) and by Vonessen
(in arbitrary characteristic), answering a question of Dixmier.
\end{abstract}

\maketitle




\section*{Introduction}


\subsection{} \label{SS:1}
Rational ideals have been rather thoroughly explored in various
settings. In the simplest case, that of an affine commutative
algebra $R$ over an algebraically closed base field $\k$, rational
ideals of $R$ are the same as maximal ideals. More generally, this
holds for any affine $\k$-algebra satisfying a polynomial identity
\cite{cP73}. For other classes of noncommutative algebras $R$,
rational ideals are identical with primitive ideals, that is,
annihilators of irreducible $R$-modules. Examples of such algebras
include group algebras of polycyclic-by-finite groups over an
algebraically closed base field $\k$ containing a non-root of unity
\cite{mLdP79} and enveloping algebras of finite-dimensional Lie
algebras over an algebraically closed field $\k$ of characteristic
$0$ \cite{cM80}, \cite{rIlS80}. Rational ideals of enveloping
algebras have been the object of intense investigation by Dixmier,
Joseph and many others from the late 1960s through the 80s; see
\S\ref{SS:history} below. The fundamental results concerning
algebraic group actions on rational ideal spectra, originally
developed in the context of enveloping algebras, were later extended
to general noetherian (or Goldie) algebras by M{\oe}glin and
Rentschler \cite{cMrR81}, \cite{cMrR84}, \cite{cMrR86},
\cite{cMrRxx} (for characteristic $0$) and by Vonessen \cite{nVE96},
\cite{nVE98} (for arbitrary characteristic). Currently, the
description of rational ideal spectra in algebraic quantum groups is
a thriving research topic; see the monograph \cite{kBkG02} by Brown
and Goodearl for an introduction. Again, rational ideals turn out to
coincide with primitive ideals for numerous examples of quantum
groups \cite[II.8.5]{kBkG02}.

\subsection{}
The aim of the present article is to liberate the theory of rational
ideals of the standard finiteness conditions, noetherianness or the
Goldie property, that are traditionally assumed in the literature.
Thus, rational ideals are defined and explored here for an arbitrary
associative algebra $R$ (with $1$) over some base field $\k$. The
Amitsur-Martindale ring of quotients will play the role of the
classical ring of quotients which underlies the usual definition of
rational ideals but need not exist in general.

Specifically, for any prime ideal $P$ of $R$, the center of the
Amitsur-Martindale ring of quotients of $R/P$, denoted by $\cC(R/P)$
and called the \emph{extended centroid} of $R/P$, is an extension
field of $\k$. The prime $P$ will be called \emph{rational} if
$$
\cC(R/P) = \k \ .
$$ 
In the special case where $R/P$ is right Goldie, $\cC(R/P)$
coincides with the center of the classical ring of quotients of
$R/P$; so our notion of rationality reduces to the familiar one in
this case. Following common practise, we will denote the collection
of all rational ideals of $R$ by $\Rat R$; so
\begin{equation*} 
\Rat R \subseteq \Spec R \ ,
\end{equation*}
where $\Spec R$ is the collection of all prime ideals of $R$, as
usual.

\subsection{}
Besides always being available, the extended centroid turns out to
lend itself rather nicely to our investigations. In fact, some of
our arguments appear to be more straightforward than earlier proofs
in more restrictive settings which were occasionally encumbered by
the fractional calculus in classical rings of quotients and by the
necessity to ensure the transfer of the Goldie property under
various constructions. Section~\ref{S:C} is preliminary in nature
and serves to deploy the definition and basic properties of extended
centroids in a form suitable for our purposes. In particular, we
show that all primitive ideals are rational under fairly general
circumstances; see Proposition~\ref{P:nullstellensatz}.

After sending out the first version of this article, we learned that
much of the material in this section was previously known, partly
even for nonassociative rings. For the convenience of the reader, we
have opted to leave our proofs intact while also indicating, to the
best of our knowledge, the original source of each result.


\subsection{}
In Section~\ref{S:Group}, we consider actions of a group $G$ by
$\k$-algebra automorphisms on $R$. Such an action induces
$G$-actions on the extended centroid $\cC(R)$ and on the set of
ideals of $R$.
Recall that a proper $G$-stable ideal $I$ of $R$ is said to be
\emph{$G$-prime} if $AB \subseteq I$ for $G$-stable ideals $A$ and
$B$ of $R$ implies that $A \subseteq I$ or $B \subseteq I$. In this
case, the subring $\cC(R/I)^G$ of $G$-invariants in $\cC(R/I)$ is an
extension field of $\k$. The $G$-prime $I$ is called
\emph{$G$-rational} if
$$
\cC(R/I)^G = \k\ .
$$
We will denote the collections of all $G$-prime and all $G$-rational
ideals of $R$ by $\GSpec R$ and $\GRat R$, respectively; so
\begin{equation*} 
\GRat R \subseteq \GSpec R \ .
\end{equation*}

The action of $G$ on the set of ideals of $R$ preserves both $\Spec
R$ and $\Rat R$. Writing the corresponding sets of $G$-orbits as $G
\backslash \Spec R$ and $G \backslash \Rat R$, the assignment $P
\mapsto \bigcap_{g \in G} g.P$ always yields a map
\begin{equation} \label{E:specmap}
G \backslash \Spec R \tto \GSpec R \ .
\end{equation}
Under fairly mild hypotheses, \eqref{E:specmap} is surjective: this
certainly holds whenever every $G$-orbit in $R$ generates a finitely
generated ideal of $R$; see Proposition~\ref{P:spec}(b). In
Proposition~\ref{P:GRat} we show that \eqref{E:specmap} always
restricts to a map
\begin{equation} \label{E:ratmap}
G \backslash \Rat R \tto \GRat R \ .
\end{equation}
More stringent conditions are required for \eqref{E:ratmap} to be
surjective. If the group $G$ is finite then \eqref{E:specmap} is
easily seen to be a bijection, and it follows from
Lemma~\ref{L:orbC} that \eqref{E:ratmap} is bijective as well.


\subsection{} \label{SS:goal}
Section~\ref{S:algebraic} focuses on rational actions of an affine
algebraic $\k$-group $G$ on $R$; the basic definitions will be
recalled at the beginning of the section. Working over an
algebraically closed base field $\k$, we show that \eqref{E:ratmap}
is then a bijection:

\begin{thm} \label{T:goal}
Let $R$ be an associative algebra over the algebraically closed
field $\k$ and let $G$ be an affine algebraic group over $\k$ acting
rationally by $\k$-algebra automorphisms on $R$. Then the map $P
\mapsto \bigcap_{g \in G} g.P$ yields a surjection $\Rat R \onto
\GRat R$ whose fibres are the $G$-orbits in $\Rat R$.
\end{thm}

The theorem quickly reduces to the situation where $G$ is connected.
Theorem~\ref{T:fibres} gives a description of the fibre of the map
$\Rat R \to \GRat R$ over any given $G$-rational ideal of $R$ for
connected $G$. This description allows us to prove transitivity of
the $G$-action on the fibres by simply invoking an earlier result of
Vonessen \cite[Theorem 4.7]{nVE98} on subfields of the rational
function field $\k(G)$ that are stable under the regular $G$-action.
Under suitable Goldie hypotheses, Theorem~\ref{T:goal} is due to
M{\oe}glin and Rentschler \cite[Th{\'e}or{\`e}me 2]{cMrRxx} in
characteristic $0$ and to Vonessen \cite[Theorem 2.10]{nVE98} in
arbitrary characteristic. The basic outline of our proof of
Theorem~\ref{T:goal} via the description of the fibres as in
Theorem~\ref{T:fibres} is adapted from the groundbreaking work of
M{\oe}glin, Rentschler and Vonessen. However, the generality of our
setting necessitates a complete reworking of the material and our
presentation contains numerous simplifications over the original
arguments.


\subsection{} \label{SS:history}
The systematics investigation of rational ideals in the enveloping
algebra $U(\fg)$ of a finite-dimensional Lie algebra $\fg$ over an
algebraically closed field $\k$ of characteristic $0$ was initiated
by Gabriel \cite{yNpG67}, \cite{pG71}. As mentioned in \S\ref{SS:1},
it was eventually established that ``rational'' is tantamount to
``primitive'' for ideals of $U(\fg)$; over an uncountable base field
$\k$, this is due to Dixmier \cite{jD77}. The reader is referred to
the standard reference \cite{jD96} for a detailed account of the
theory of primitive ideals in enveloping algebras; for an updated
survey, see \cite{rR87}. Here we just mention that the original
motivation behind Theorem~\ref{T:goal} and its predecessors was a
question of Dixmier \cite{jD72} (see also \cite[Problem 11]{jD96})
concerning primitive ideals of $U(\fg)$. Specifically, if $G$ is the
adjoint algebraic group of $\fg$ then, for any ideal $\fk$ of $\fg$
and any primitive ideal $Q$ of $U(\fg)$, the ideal $I = Q \cap
U(\fk)$ of $U(\fk)$ is $G$-rational \cite{jD77}. Dixmier asked if
the following are true for $I$:
\begin{enumerate}
\item $I = \bigcap_{g \in G} g.P$ for some primitive ideal
$P$ of $U(\fk)$, and
\item any two such primitive ideals belong to the same $G$-orbit.
\end{enumerate}
The earlier version of Theorem~\ref{T:goal}, due to M{\oe}glin and
Rentsch\-ler, settled both (a) and (b) in the affirmative. Letting
$\Prim U(\fk)$ denote the collection of all primitive ideals of
$U(\fk)$ endowed with the Jacobson-Zariski topology, (a) says that
the set $\{J \in \Prim U(\fk) \mid J \supseteq I \}$ is the closure
of the orbit $G.P$ in $\Prim U(\fk)$. Following Dixmier \cite{jD72}
such $P$ are called \emph{generic} for $I$. The uniqueness of
generic orbits as in (b) was proved for solvable $\fg$ in
\cite{BGR73} and generally (over uncountable $\k$) in \cite{rR79};
this fact was instrumental for the proof that the Dixmier and Duflo
maps are injective in the solvable and algebraic case, respectively
(Rentschler \cite{rR74}, Duflo \cite{mD82}).


\subsection{} \label{SS:future}
In future work, we hope to address some topological aspects of $\Rat
R$ endowed with the Jacobson-Zariski topology from $\Spec R$.
Finally, it remains to bring the machinery developed herein to bear
on new classes of algebras that lack the traditional finiteness
conditions.




\section{The extended centroid} \label{S:C}

Throughout this section, $R$ will denote an associative ring. It is
understood that all rings have a $1$ which is inherited by subrings
and preserved under homomorphisms.


\subsection{The Amitsur-Martindale ring of quotients} \label{SS:Qrt}
Let $\sE = \sE(R)$ denote the filter consisting of all (two-sided)
ideals $I$ of $R$ such that
$$
\lann_RI = \{ r \in R \mid rI = 0 \} = 0 \ .
$$
The right Amitsur-Martindale ring of quotients, introduced for prime
rings $R$ by Martidale \cite{wM69} and in general by Amitsur
\cite{sA72}, is defined by
$$
\Qr(R) = \varinjlim_{I \in \sE} \Hom(I_R,R_R) \ .
$$
Explicitly, the elements of $\Qr(R)$ are equivalence classes of
right $R$-module maps $f \colon I_R \to R_R$ with $I \in \sE$; the
map $f$ is defined to be equivalent to $f' \colon I'_R \to R_R$ $(I'
\in \sE)$ if $f$ and $f'$ agree on some ideal $J \subseteq I \cap
I'$, $J \in \sE$. In this case, $f$ and $f'$ actually agree on $I
\cap I'$; see \cite[Lemma 1]{sA72}. The sum of two elements $q,q'
\in \Qr(R)$, represented by $f \colon I_R \to R_R$ $(I \in \sE)$ and
$f' \colon I'_R \to R_R$ $(I' \in \sE)$, respectively, is defined to
be the class of $f + f' \colon I \cap I' \to R$. Similarly, the
product $qq' \in \Qr(R)$ is the class of the composite $f \circ f'
\colon I'I \to R$. This makes $\Qr(R)$ into a ring; the identity
element is the class of the identity map $\Id_R$ on $R$. Sending an
element $r \in R$ to the equivalence class of the map $\lambda_r
\colon R \to R$, $x \mapsto rx$, yields an embedding of $R$ as a
subring of $\Qr(R)$. Suppose the element $q \in \Qr(R)$ is
represented by $f \colon I_R \to R_R$ $(I \in \sE)$. Then the
equality $f \circ \lambda_r = \lambda_{f(r)}$ $(r \in I)$ shows that
$q I \subseteq R$.

We summarize the foregoing and some easy consequences thereof in the
following proposition. Complete details can be found in \cite{sA72}
and in \cite[Proposition 10.2]{dP89}, for example.

\begin{prop}
\label{P:Qr}
The ring $\Qr(R)$ has the following properties:
\begin{enumerate}
\item[i.] There is a ring embedding $R \into \Qr(R)$;
\item[ii.] for each $q \in \Qr(R)$, there exits $I \in \sE$ with
$qI \subseteq R$;
\item[iii.] if $qI = 0$ for $q \in \Qr(R)$ and $I \in \sE$ then
$q = 0$;
\item[iv.] given $f \colon I_R \to R_R$ with $I \in \sE$, there exists
$q \in \Qr(R)$ with $qr = f(r)$ for all $r \in I$.
\end{enumerate}
Furthermore, $\Qr(R)$ is characterized by these properties: any
other ring satisfying {\rm (i)} -- {\rm (iv)} is $R$-isomorphic to
$\Qr(R)$.
\end{prop}


\subsection{The extended centroid} \label{SS:C}

The extended centroid of $R$ is defined to be the center of
$\Qr(R)$; it will be denoted by $\cC(R)$:
$$
\cC(R) = \cen(\Qr(R)) \ .
$$
It is easy to see from Proposition~\ref{P:Qr} that $\cC(R)$
coincides with the centralizer of $R$ in $\Qr(R)$:
$$
\cC(R) = C_{\Qr(R)}(R) = \{ q \in \Qr(R) \mid qr=rq \ \forall r \in
R \} \ .
$$
In particular, the center $\cen(R)$ of $R$ is contained in $\cC(R)$.
Moreover, an element $q \in \Qr(R)$ belongs to $\cC(R)$ if and only
if $q$ is represented by an $(R,R)$-\emph{bimodule} map $f \colon I
\to R$ with $I \in \sE$; in this case, every representative $f'
\colon I'_R \to R_R$ $(I' \in \sE)$ of $q$ is an $(R,R)$-bimodule
map; see \cite[Theorem 3]{sA72}.

\subsubsection{} \label{SSS:leftright}
By reversing sides, one can define the left ring of quotients
$\Ql(R)$ and its center $\cCl(R) = \cen(\Ql(R))$ as above. However,
we will mainly be concerned with semiprime rings, that is, rings $R$
having no nonzero ideals of square $0$. In that case, $\lann_RI =
\rann_RI$ holds for every ideal $I$ of $R$; so the definition  of
$\sE(R)$ is symmetric. Moreover, any $q \in \cC(R)$ is represented
by an $(R,R)$-bimodule map $f \colon I \to R$ with $I \in \sE$. The
class of $f$ in $\Ql(R)$ is an element $q' \in \cCl(R)$, and the map
$q \mapsto q'$ yields an isomorphism $\cC(R) \iso \cCl(R)$. In the
following, we shall always work with $\Qr(R)$ and $\cC(R)$.

\subsubsection{}
Let $R$ be semiprime. Then one knows that $\cC(R)$ is a von Neumann
regular ring. Moreover, $R$ is prime if and only if $\cC(R)$ is a
field; see \cite[Theorem 5]{sA72}.


\subsection{Central closure} \label{SS:Central}

Rings $R$ such that $\cC(R) \subseteq R$ are called \emph{centrally
closed}. In this case, $\cC(R) = \cen(R)$. For every semiprime ring
$R$, the subring $R\cC(R)$ of $\Qr(R)$ is a semiprime centrally
closed ring called the \emph{central closure} of $R$; see
\cite[Theorem 3.2]{wBwM79}. If $R$ is prime then so is the central
closure $R\cC(R)$ by Proposition~\ref{P:Qr}(ii).

The following lemma goes back to Martindale \cite{wM69}.

\begin{lem} \label{L:Central}
Let $R$ be a prime centrally closed ring and let $S$ be an algebra
over the field $C = \cC(R)$. Then:
\begin{enumerate}
\item
Every nonzero ideal $I$ of $R \otimes_C S$ contains an element $0
\neq r \otimes s$ with $r \in R$, $s \in S$.
\item
If $S$ is simple then every nonzero ideal $I$ of $R \otimes_C S$
intersects $R$ nontrivially. Consequently, $R \otimes_C S$ is prime.
\item
If $I$ is a prime ideal of $R \otimes_C S$ such that $I \cap R = 0$
then $I = R \otimes_C (I \cap S)$.
\end{enumerate}
\end{lem}

\begin{proof}
(a)
Fix a $C$-basis $\{ s_i \}$ of $S$. Consider an element $0 \neq t
= \sum_i r_i \otimes s_i \in I$ with a minimal number of nonzero
$R$-coefficients $r_i$ among all nonzero elements of $I$ and choose
$i_0$ with $r = r_{i_0} \neq 0$. Then the element $rxt - txr =
\sum_{i\neq i_0} (rxr_i - r_ixr) \otimes s_i$ must be zero for all
$x \in R$. Hence $rxr_i = r_ixr$ holds for all $i$, and by
\cite[Theorem 1]{wM69}, there are $c_i \in C$ such that $r_i =
rc_i$. Therefore, $t = r \otimes s$ with $s = \sum_i c_i s_i \in S$.

(b) If $S$ is simple then we can make $s=1$ in (a), and so $0 \neq r
\in I \cap R$. Since $R$ is prime, it follows that $R \otimes_C S$
is prime as well.

(c) Suppose for a contradiction that $I \supsetneq R \otimes_C (I
\cap S)$. Replacing $S$ by $S/(I \cap S)$, we may assume that $I
\neq 0$ but $I \cap R = 0$ and $I \cap S = 0$. Choosing $r \otimes s
\in I$ as in (a), we obtain that $I \supseteq S(r \otimes s)R = rR
\otimes_C Ss$. Since $I$ is prime, we must have $r \in I$ or $s \in
I$, whence the desired contradiction.
\end{proof}


\subsection{Examples} \label{SS:Ex}

\subsubsection{} \label{SSS:ex1}
If $R$ is a simple ring, or a finite product of simple rings, then
$\sE(R) = \{R\}$, and hence $\Qr(R) = R$ by
Proposition~\ref{P:Qr}(i)(ii). Thus, $R$ is certainly centrally
closed in this case. Less trivial examples of centrally closed rings
include crossed products $R*F$ with $R$ a simple ring and $F$ a free
semigroup on at least two generators (\cite[Theorem 13.4]{dP89}) and
Laurent power series rings $R((x))$ over centrally closed rings $R$
(\cite{wMmRjR90}).

\subsubsection{} \label{SSS:semiprimeGoldie}
If $R$ is semiprime right Goldie then $\cC(R) = \cen(\Qcl(R))$, the
center of the classical ring of quotients of $R$. Indeed, $\Qcl(R)$
coincides with the maximal ring of quotients $\Qmax(R)$ in this
case; see, e.g., Lambek \cite[Prop.~4.6.2]{jL76}. Furthermore, the
Amitsur-Martindale ring of quotients $\Qr(R)$ is $R$-isomorphic to
the subring of $\Qmax(R)$ consisting of all $q \in \Qmax(R)$ such
that $q I \subseteq R$ for some $I \in \sE(R)$; see, e.g.,
\cite[Chap.~24]{dP91} or \cite[Chap.~3]{sM80}. This isomorphism
yields an isomorphism $\cC(R) \cong \cen(\Qmax(R))$.

\subsubsection{} \label{SSS:Ug}
Let $R$ be a semiprime homomorphic image of the enveloping algebra
$U(\fg)$ of a finite-dimensional Lie algebra $\fg$ over some base
field $\k$.
Answering a question of Rentschler, we show here that
\begin{equation*}
\text{\emph{$\Qr(R)$ consists of all $\ad\fg$-finite elements of
$\Qcl(R)$.}}
\end{equation*}
Here, $\ad \colon U(\fg) \to \End_\k \Qcl(R)$ is the standard
adjoint action, given by $\ad x(q) = xq - qx$ for $x \in \fg$ and $q
\in \Qcl(R)$, and $q$ is called $\ad\fg$-finite if the $\k$-subspace
$\ad U(\fg)(q)$ of $\Qcl(R)$ is a finite dimensional.
To prove the claim, recall from \S\ref{SSS:semiprimeGoldie} that
$\Qr(R) = \{ q \in \Qcl(R) \mid \text{$q I \subseteq R$ for some $I
\in \sE(R)$} \}$. First consider $q \in \Qr(R)$. Letting $R_n$ and
$I_n = I \cap R_n$ $(n\ge 0)$ denote the filtrations of $R$ and $I$,
respectively, that are induced by the canonical filtration of
$U(\fg)$ (\cite[2.3.1]{jD96}), we have $I = I_sR$ and $qI_s
\subseteq R_t$ for suitable $s,t \ge 0$. Since both $I_s$ and $R_t$
are $\ad(\fg)$-stable, it follows that $\ad U(\fg)(q) I_s \subseteq
R_t$. Furthermore, $\lann_{\Qcl(R)}I_s = \lann_{\Qcl(R)}I = 0$; so
$\ad U(\fg)(q)$ embeds into $\Hom_{\k}(I_s,R_t)$ proving that $q$ is
$\ad\fg$-finite. Conversely, suppose that $q \in \Qcl(R)$ is
$\ad\fg$-finite and let $\{ q_i\}_1^m$ be a $\k$-basis of $\ad
U(\fg)(q)$. Each $D_i = \{ r \in R \mid q_ir \in R \}$ is an
essential right ideal of $R$, and hence $I = \bigcap_1^m D_i = \{ r
\in R \mid \ad U(\fg)(q) r \subseteq R \}$ is an essential right
ideal of $R$ which is also $\ad(\fg)$-stable, since this holds for
$\ad U(\fg)(q)$ and $R$. Therefore, $I \in \sE(R)$ which shows that
$q \in \Qr(R)$.


\subsection{Centralizing homomorphisms} \label{SS:Hom}

A ring homomorphism $\phi \colon R \to S$ is called
\emph{centralizing} if the ring $S$ is generated by $\phi(R)$ and
the centralizer $C_S(\phi(R)) = \{ s \in S \mid s \phi(r) = \phi(r)s
\ \forall r \in R \}$. Surjective ring homomorphisms are clearly
centralizing, and composites of centralizing homomorphisms are again
centralizing. Note also that any centralizing homomorphism $\phi
\colon R \to S$ sends the center $\cen(R)$ of $R$ to $\cen(S)$.
Finally, $\phi$ induces a map $\Spec S \to \Spec R$, $P \mapsto
\phi^{-1}(P)$.

For any $q \in \Qr(R)$, we define the ideal $D_q$ of $R$ by
\begin{equation} \label{E:Dq}
D_q = \{ r \in R \mid q R r \subseteq R \} \ .
\end{equation}
By Proposition~\ref{P:Qr}(ii), $D_q \in \sE(R)$. If $q \in \cC(R)$
then the description of the ideal $D_q$ simplifies to $D_q = \{ r
\in R \mid q r \subseteq R \}$.

\begin{lem} \label{L:C}
Let $\phi \colon R \to S$ be a centralizing homomorphism of rings.
Put
$$
\cC_\phi = \{ q \in \cC(R) \mid \lann_S\phi(D_q) = 0 \} \ .
$$
Then $R\cC_\phi$ is a subring of $\Qr(R)$ containing $R$. The map
$\phi$ extends uniquely to a centralizing ring homomorphism
$\til{\phi} \colon R\cC_\phi \to S\cC(S)$. In particular,
$\til{\phi}(\cC_\phi) \subseteq \cC(S)$.
\end{lem}

\begin{proof}
Put
$$
R_\phi = \{ q \in \Qr(R) \mid \lann_S\phi(D_q) = 0 \}\ .
$$
Since $R = \{ q \in \Qr(R) \mid 1 \in D_q \}$, we certainly have $R
\subseteq R_\phi$. For  $q,q' \in \Qr(R)$, one easily checks that
$D_{q'}D_q \subseteq D_q \cap D_{q'} \subseteq D_{q+q'}$ and
$D_{q'}D_q \subseteq D_{qq'}$. Moreover, if $\phi(D_q)$ and
$\phi(D_{q'})$ both have zero left annihilator in $S$ then so does
$\phi(D_{q'}D_q ) = \phi(D_{q'})\phi(D_q)$. This shows that $q + q'
\in R_\phi$ and $q q' \in R_\phi$ for $q,q' \in R_\phi$; so $R_\phi$
is a subring of $\Qr(R)$ containing $R$. Since $\cC_\phi =
\cen(R_\phi)$, it follows that $R\cC_\phi$ is also a subring of
$\Qr(R)$ containing $R$.

Now let $q \in \cC_\phi$ be given. Then $\phi(D_q)S =
\phi(D_q)C_S(\phi(R)) \in \sE(S)$. Define $\bar{q} \colon \phi(D_q)S
\to S$ by
\begin{equation*}
\bar{q}(\sum_i \phi(x_i)c_i) = \sum_i \phi(qx_i)c_i
\end{equation*}
for $x_i \in D_q$, $c_i \in C_S(\phi(R))$. To see that $\bar{q}$ is
well-defined, note that, for each $d \in D_q$, we have
\begin{equation*}
\begin{split}
\sum_i \phi(x_i)c_i \phi(qd) &= \sum_i \phi(x_i)\phi(qd) c_i =
\sum_i \phi(x_i q d) c_i \\
&= \sum_i \phi(q x_i d) c_i = \sum_i \phi(q x_i)\phi(d) c_i \\
&= \sum_i
\phi(q x_i)c_i \phi(d) \ .
\end{split}
\end{equation*}
Thus, if $\sum_i \phi(x_i)c_i = \sum_j \phi(y_j)e_j$ with $x_i,y_j
\in D_q$ and $c_i,e_j \in C_S(\phi(R))$ then the above computation
gives
$$
0 = \left( \sum_i \phi(x_i)c_i - \sum_j \phi(y_j)e_j \right) \phi(q
D_q) = \left( \sum_i \phi(qx_i)c_i - \sum_j \phi(qy_j)e_j \right)
\phi(D_q) \ ,
$$
and so $0 = \sum_i \phi(qx_i)c_i - \sum_j \phi(qy_j)e_j$. Therefore,
$\bar{q}$ is well-defined.

It is straightforward to check that $\bar{q}$ is an $(S,S)$-bimodule
map. Hence, the class of $\bar{q}$ in $\Qr(R)$ is an element
$\til{\phi}(q) \in \cC(S)$. The map $q \mapsto \til{\phi}(q)$ is a
ring homomorphism $\cC_\phi \to \cC(S)$ which yields the desired
extension $\til{\phi} \colon R\cC_\phi \to S\cC(S)$:
$\til{\phi}(\sum_i r_i q_i) = \sum_i \phi(r_i)\til{\phi}(q_i)$ for
$r_i \in R$, $q_i \in \cC_\phi$. Well-definedness and uniqueness of
$\til{\phi}$ follow easily from the fact that, given finitely many
$x_i \in R_\phi$, there is an ideal $D$ of $R$ with $\lann_S\phi(D)
= 0$ and $x_iD \subseteq R$ for all $i$.
\end{proof}

In the special case where both $R$ and $S$ are commutative domains
in Lemma~\ref{L:C} above, we have $\Qr(R) = \cC(R) = \Fract R$, the
classical field of fractions of $R$, and similarly for $S$.
Moreover, $R\cC_\phi = R_P$ is the localization of $R$ at the prime
$P = \Ker \phi$ and the map $R\cC_\phi \to S\cC(S)$ is the usual map
$R_P \to \Fract S$.


\subsection{Extended centroids and primitive ideals} \label{SS:primitive}

By Schur's Lemma, the endomorphism ring $\End_RV$ of any simple
$R$-module $V_R$ is a division ring. The following lemma is
well-known in the special case of noetherian (or Goldie) rings (see,
e.g., Dixmier \cite[4.1.6]{jD96}); for general rings, the lemma was
apparently first observed by Martindale \cite[Theorem 12]{wm69b}.
Since the latter result is stated in terms of the so-called complete
ring of quotients, we include the proof for the reader's
convenience.

\begin{lem} \label{L:primitive}
Let $V_R$ be a simple $R$-module, and let $P = \ann_RV$ be its
annihilator. Then the canonical embedding $\cen(R/P) \into
\cen\left( \End_RV \right)$ extends to an embedding of fields
$$
\cC(R/P) \into \cen\left( \End_RV \right) \ .
$$
\end{lem}

\begin{proof}
We may assume that $P=0$. For a given $q \in \cC(R)$, we wish to
define an endomorphism $\delta_q \in \cen\left( \End_RV \right)$. To
this end, note that every $x \in V$ can be written as $x = vd$ for
suitable $d \in D_q$, $v \in V$. Define
\begin{equation*} 
\delta_q(x) = v(dq) \in V \ .
\end{equation*}
To see that this is well-defined, assume that $vd = v'd'$ holds for
$v,v' \in V$ and $d,d' \in D_q$. Then $\left( v(dq) - v'(d'q)
\right)D_q = \left( vd - v'd' \right)(qD_q) = 0$ and so $v(dq) -
v'(d'q) = 0$. It is straightforward to check that $\delta_q \in
\End_RV$. Moreover, for any $\delta \in \End_RV$ and $vd \in V$, one
computes
$$
\delta\delta_q(vd) = \delta(v(dq)) =  \delta(v)(dq) =
\delta_q(\delta(v)d) = \delta_q\delta(vd)  \ .
$$
Thus, $\delta_q \in \cen\left( \End_RV \right)$. The map $\cC(R) \to
\cen\left( \End_RV \right)$, $q \mapsto \delta_q$, is easily seen to
be additive. Furthermore, for $q,q' \in \cC(R)$, $d \in D_q$, $d'
\in D_{q'}$ and $v \in V$, one has
$$
\delta_{qq'}(vd'd) = v(d'dqq') = v(d'q')(dq) =
\delta_{q}(\delta_{q'}(vd')d) = \delta_{q}(\delta_{q'}(vd'd)) \ .
$$
Thus, the map is a ring homomorphism; it is injective because
$\cC(R)$ is a field.
\end{proof}


\subsection{Rational algebras and ideals} \label{SS:rat}

An algebra $R$ over some field $\k$ will be called \emph{rational}
(or \emph{$\k$-rational}) if $R$ is prime and $\cC(R) = \k$. A prime
ideal $P$ of $R$ will be called rational if $R/P$ is a rational
$\k$-algebra. In view of \S\ref{SSS:leftright}, the notion of
rationality is left-right symmetric.

We remark that rational $\k$-algebras are called \emph{closed over
$\k$} in \cite{EMO75} where such algebras are investigated in a
non-associative context. Alternatively, one could define a prime
$\k$-algebra $R$ to be rational if the field extension $\cC(R)/\k$
is algebraic; for noetherian (or Goldie) algebras, this version of
rationality is adopted in many places in the literature (e.g.,
\cite{kBkG02}). However, we will work with the above definition
throughout.

\subsubsection{} \label{SSS:rat1}
By \S\ref{SS:Central} the central closure $R\cC(R)$ of any prime
ring $R$ is $\cC(R)$-rational.

\subsubsection{} \label{SSS:rat2}
The Schur division rings $\End_RV$ considered in
\S\ref{SS:primitive} are division algebras over $\k$, and their
centers are extension fields of $\k$. We will say that the algebra
$R$ satisfies the \emph{weak Nullstellensatz} if $\cen\left( \End_RV
\right)$ is algebraic over $\k$ for every simple $R$-module $V_R$.

\begin{prop} \label{P:nullstellensatz}
If $R$ is a $\k$-algebra satisfying the weak Nullstellensatz and
$\k$ is agebraically closed then all primitive ideals of $R$ are
rational.
\end{prop}

\begin{proof}
By hypothesis, $\cen\left( \End_RV \right) = \k$ holds for every
simple $R$-module $V_R$. It follows from Lemma~\ref{L:primitive}
that $P = \ann_RV$ satisfies $\cC(R/P) = \k$.
\end{proof}

For an affine commutative $\k$-algebra $R$, the Schur division
algebras in question are just the quotients $R/P$, where $P$ is a
maximal ideal of $R$. The classical weak Nullstellensatz is
equivalent to the statement that $R/P$ is always algebraic over
$\k$; see, e.g., Lang \cite[Theorem IX.1.4]{sL02}. Thus affine
commutative algebras do satisfy the weak Nullstellensatz.

Many noncommutative algebras satisfying the weak Nullstellensatz are
known; see \cite[Chapter 9]{jMcCjR87} for an overview. In fact, as
long as the cardinality of the base field $\k$ is larger than
$\dim_{\k}R$, the weak Nullstellensatz is guaranteed to hold; see
\cite[Corollary 9.1.8]{jMcCjR87} or \cite[II.7.16]{kBkG02}. This
applies, for example, to any countably generated algebra over an
uncountable field $\k$.


\subsection{Scalar extensions} \label{SS:scalar}

We continue to let $R$ denote an algebra over some field $\k$. For
any given $\k$-algebra $A$, we have an embedding
$$
\Qr(R) \otimes_\k A \into \Qr(R \otimes_\k A)
$$
which extends the canonical embedding $R \otimes_\k A \into \Qr(R
\otimes_\k A)$. For, let $q \in \Qr(R)$ be represented by the map $f
\colon I_R \to R_R$ with $I \in \sE(R)$. Then $I \otimes_\k A \in
\sE(R \otimes_\k A)$. Sending $q$ to the class of the map $f \otimes
\Id_A$ we obtain a ring homomorphism $\Qr(R) \to \Qr(R \otimes_\k
A)$ extending the canonical embedding $R \into R \otimes_\k A \into
\Qr(R \otimes_\k A)$. By Proposition~\ref{P:Qr}(ii),(iii), the image
of $\Qr(R)$ in $\Qr(R \otimes_\k A)$ commutes with $A$ and the
resulting map $\Qr(R) \otimes_\k A \to \Qr(R \otimes_\k A)$ is
injective. Moreover, since $f \otimes \Id_A$ is an $(R \otimes_\k
A,R \otimes_\k A)$-bimodule map if $f$ is an $(R,R)$-bimodule map,
the embedding of $\Qr(R)$ into $\Qr(R \otimes_\k A)$ sends $\cC(R)$
to $\cC(R \otimes_\k A)$. Thus, if $A$ is commutative, this yields
an embedding
\begin{equation} \label{E:scalar}
\cC(R) \otimes_\k A \into \cC(R \otimes_\k A)\ .
\end{equation}

The following lemma is the associative case of \cite[Theorem
3.5]{EMO75}.

\begin{lem} \label{L:scalar}
Assume that $R$ is rational. Then, for every field extension $K/\k$,
the $K$-algebra $R_K = R \otimes_\k K$ is rational.
\end{lem}

\begin{proof}
By Lemma~\ref{L:Central}(b), we know that $R_K$ is prime. Moreover,
for any given $q \in \cC(R_K)$, we may choose an element $0 \neq x
\in D_q \cap R$. Fix a $\k$-basis $\{ k_i \}$ for $K$. The map
$$
q_i \colon I = RxR \stackrel{q \cdot}{\tto} R_K
\stackrel{\text{proj}}{\onto} R \otimes k_i \iso R
$$
is an $(R,R)$-bimodule map. Hence $q_i$ is multiplication with some
$c_i \in \k$, by hypothesis on $R$, and all but finitely many $c_i$
are zero. Therefore, the map $I \stackrel{q \cdot}{\tto} R_K$ is
multiplication with $k = \sum_i c_i k_i \in K$. Consequently, $q = k
\in K$.
\end{proof}


\section{Group actions} \label{S:Group}

In this section, we assume that a group $G$ acts by automorphisms on
the ring $R$; the action will be written as $G \times R \to R$,
$(g,r) \mapsto g.r$.


\subsection{} \label{SS:orb}
Let $M$ be a set with a left $G$-action $G \times M \to M$, $(g,m)
\mapsto g.m$. For any subset $X$ of $M$,
$$
G_X = \stab_GX = \{ g \in G \mid g.X = X \}
$$
will denote the isotropy group of $X$. Furthermore, we put
$$
(X:G) = \bigcap_{g \in G} g.X \ ;
$$
this is the largest $G$-stable subset of $M$ that is contained in
$X$. 
We will be primarily concerned with the situation where $M=R$ and
$X$ is an ideal of $R$ in which case $(X:G)$
is also an ideal of $R$. 

\subsection{$G$-primes} \label{SS:GSpec}
The ring $R$ is said to be \emph{$G$-prime} if $R \neq 0$ and the
product of any two nonzero $G$-stable ideals of $R$ is again
nonzero. A $G$-stable ideal $I$ of $R$ is called $G$-prime if $R/I$
is a $G$-prime ring for the $G$-action on $R/I$ coming from the
given action of $G$ on $R$. In the special case where the $G$-action
on $R$ is trivial, $G$-primes of $R$ are just the prime ideals of
$R$ in the usual sense. Recall that the collection of all $G$-prime
ideals of $R$ is denoted by $\GSpec R$ while $\Spec R$ is the
collection of all ordinary primes of $R$.

\begin{prop} \label{P:spec}
\begin{enumerate}
\item
There is a well-defined map
\begin{equation*} 
\Spec R \tto \GSpec R \ ,\qquad P \mapsto (P:G)\ .
\end{equation*}
\item Assume that,  for each $r \in R$, the $G$-orbit $G.r$ generates a
finitely generated ideal of $R$. Then the map in (a) is surjective.
In particular, all $G$-primes of $R$ are semiprime in this case.
\end{enumerate}
\end{prop}

\begin{proof}
It is straightforward to check that $(P:G)$ is $G$-prime for any
prime ideal $P$ of $R$; so (a) is clear.

For (b), consider a $G$-prime ideal $I$ of $R$. We will show that
there is an ideal $P$ of $R$ which is maximal subject to the
condition $(P:G) = I$; the ideal $P$ is then easily seen to be
prime. In order to prove the existence of $P$, we use Zorn's Lemma.
So let $\{ I_j \}$ be a chain of ideals of $R$ such that $(I_j:G) =
I$ holds for all $j$. We need to show that the ideal $I_* =
\bigcup_j I_j$ satisfies $(I_*:G) = I$. For this, let $r \in
(I_*:G)$ be given. Then the ideal $(G.r)$ that is generated by $G.r$
is contained in $(I_*:G)$ and $(G.r)$ is a finitely generated
$G$-stable ideal of $R$. Therefore, $(G.r) \subseteq (I_{j}:G)$ for
some $j$ and so $r \in I$, as desired.
\end{proof}

For brevity, we will call $G$-actions satisfying the finiteness
hypothesis in (b) above \emph{locally ideal finite}. Clearly, all
actions of finite groups as well as all group actions on noetherian
rings are locally ideal finite. Another important class of examples
are the \emph{locally finite} actions in the usual sense: by
definition, these are $G$-actions on some $\k$-algebra $R$ such that
the $G$-orbit of each $r \in R$ generates a finite-dimensional
$\k$-subspace of $R$. This includes the rational actions of
algebraic groups to be considered in Section~\ref{S:algebraic}. In
all these cases, Proposition~\ref{P:spec} is a standard result; the
argument given above is merely a variant of earlier proofs.


\subsection{$G$-primes and the extended centroid} \label{SS:GSpecC}

The $G$-action on $R$ extends uniquely to an action of $G$ on
$\Qr(R)$: if $q \in \Qr(R)$ is represented by $f \colon I_R \to R_R$
$(I \in \sE)$ then $g.q \in \Qr(R)$ is defined to be the class of
the map $g.f \colon g.I \to R$ that is given by $(g.f)(g.x) =
g.f(x)$ for $x \in I$. Therefore, $G$ also acts on the extended
centroid $\cC(R)$ of $R$. As usual, the ring of $G$-invariants in
$\cC(R)$ will denoted by $\cC(R)^G$.

\begin{prop} \label{P:Gprime}
If $R$ is $G$-prime then $\cC(R)^G$ is a field. Conversely, if $R$
is semiprime and $\cC(R)^G$ is a field then $R$ is $G$-prime.
\end{prop}

\begin{proof} We follow the outline of the proof of \cite[Theorem 5]{sA72}.

First assume that $R$ is $G$-prime and let $0 \neq q \in \cC(R)^G$
be given. Then $qD_q$ is a nonzero $G$-stable ideal of $R$, and
hence $\lann_R(qD_q) = 0$ because $R$ is $G$-prime. So $qD_q \in
\sE(R)$. Moreover, $\ann_R(q) = \{ r \in R \mid rq = 0 \} \subseteq
\lann_R(qD_q)$ and so $\ann_R(q) = 0$. Therefore, the map $D_q \to
qD_q$, $r \mapsto qr = rq$, is an $(R,R)$-bimodule isomorphism which
is $G$-equivariant. The class of the inverse map belongs to
$\cC(R)^G$ and is the desired inverse for $q$.

Next, assume that $R$ is semiprime but not $G$-prime. Then there
exists a nonzero $G$-stable ideal $I$ of $R$ such that $J =
\lann_R(I) \neq 0$. Since $R$ is semiprime, the sum $I + J$ is
direct and $I + J \in \sE(R)$. Define maps $f,f' \colon I + J \to R$
by $f(i+j) = i$ and $f'(i+j) = j$. Letting $q$ and $q'$ denote the
classes of $f$ and $f'$, respectively, in $\Qr(R)$ we have $f,f' \in
\cC(R)^G$ and $ff' = 0$. Therefore, $\cC(R)^G$ is not a field.
\end{proof}

The following technical lemma
will be crucial. Recall that $G_I$ denotes the isotropy group of
$I$.

\begin{lem} \label{L:orbC}
Let $P$ be a prime ideal of $R$. \begin{enumerate} \item For every
subgroup of $H \le G$, the canonical map $R/(P:G) \onto R/(P:H)$
induces an embedding of fields
$$
\cC(R/(P:G))^G \into \cC(R/(P:H))^{\Htil} \ .
$$
The degree of the field extension is at most $[G : \Htil]$.
\item
If $P$ has a finite $G$-orbit then we obtain an isomorphism of
fields
\begin{equation*} 
\cC(R/(P:G))^G \iso \cC(R/P)^{G_P}\ .
\end{equation*}
\end{enumerate}
\end{lem}

\begin{proof}
(a) After factoring out the ideal $(P:G)$ we may assume that
$(P:G)=0$, $R$ is $G$-prime, and $\cC(R)^G$ is a field; see
Propositions~\ref{P:spec} and \ref{P:Gprime}.
Consider the canonical map
$$
\phi \colon R \onto S:= R/(P:H) \ .
$$
Using the notation of Lemma~\ref{L:C}, we have $\cC(R)^G \subseteq
\cC_\phi$. Indeed, for each $q \in \cC(R)^G$, the ideal $D_q$ is
nonzero and $G$-stable, and hence $D_q \nsubseteq P$. Therefore,
$\phi(D_q)$ is a nonzero $H$-stable ideal of the $H$-prime ring $S$,
and so $\phi(D_q) \in \sE(S)$. The map $\cC_\phi \to \cC(S)$
constructed in Lemma~\ref{L:C} yields an embedding embedding
$\cC(R)^G \into \cC(S)$: the image of $q \in \cC(R)^G$ is the class
of the map $f \colon \phi(D_q) \to S$ that is defined by
$f(\phi(x)) = \phi(qx)$ for $x \in D_q$. Since $\phi$ is
$\Htil$-equivariant, one computes, for $x \in D_q$ and $g \in
\Htil$,
\begin{equation*}
\begin{split}
(g.f)(g.\phi(x)) &= g.f(\phi(x)) = g.\phi(qx) = \phi(g.(qx)) \\
&= \phi(q(g.x)) = f(\phi(g.x)) = f(g.\phi(x))\ ;
\end{split}
\end{equation*}
so $g.f = f$. Therefore the image of $\cC(R)^G$ is
contained in $\cC(S)^{\Htil}$.

It remains to show that $[\cC(S)^{\Htil} : \cC(R)^G] \le [G:\Htil]$
if the latter number is finite. To this end, put $N = \bigcap_{g \in
G} g^{-1} \Htil g$; this is a normal subgroup of $G$ which has finite index
in $G$ and is contained in $\Htil$. Since $(P:H) = (P:\Htil)$, the
foregoing yields embeddings of fields
$$
\cC(R)^G \into \cC(S)^{\Htil} = \cC(R/(P:\Htil))^{\Htil} \into
\cC(R/(P:N))^{N'}\ ,
$$
where $N':= \Htil \cap G_{(P:N)}$. The image of $\cC(R)^G$ under the
composite embedding is contained in $\cC(R/(P:N))^{G_{(P:N)}}$ and,
by Galois theory,
$$
[\cC(R/(P:N))^{N'} : \cC(R/(P:N))^{G_{(P:N)}}] \le
[G_{(P:N)} :N'] \le [G : \Htil]\ .
$$
It suffices to show that the image of $\cC(R)^G$ is actually
equal to $\cC(R/(P:N))^{G_{(P:N)}}$. Therefore, replacing $H$ by
$N$, it suffices to show:
\begin{equation} \label{E:orbC}
\text{If $H \unlhd G$ and $[G:\Htil] < \infty$ then $\cC(R)^G$ maps
onto $\cC(S)^{\Htil}$.}
\end{equation}
To this end, we will prove the following
\begin{claim}
Let $t \in \cC(S)^{\Htil}$ be given. There exists a $G$-stable ideal
$I$ of $R$ such that $0 \neq \phi(I) \subseteq D_t$ and such that,
for every $x \in I$, there exists an $x' \in R$ satisfying
\begin{equation} \label{E:GSpecC}
\phi(g.x') = t\phi(g.x) \quad \text{for all $g \in G$.}
\end{equation}
\end{claim}
Note that $G$-stability of $I$ and the condition $\phi(I) \subseteq
D_t$ ensure that $t\phi(g.x) \in S$ holds for all $g \in G$, $x \in
I$. Moreover, any $G$-stable ideal $I$ satisfying $0 \neq \phi(I)$
belongs to $\sE(R)$. For, $\lann_S\phi(I) = 0$ since $S$ is
$H$-prime, and hence $\lann_RI$ is contained in $(P:G) = 0$.
Finally,  the element $x'$ is uniquely determined by
\eqref{E:GSpecC} for any given $x$. Indeed, if $x'' \in R$ also
satisfies \eqref{E:GSpecC} then $\phi(g.x') = \phi(g.x'')$ holds for
all $g \in G$ and so $x' - x'' \in (P:G) = 0$. Therefore, assuming
the claim for now, we can define a map
$$
f \colon I \to R\ ,\quad x \mapsto x' \ .
$$
It is easy to check that $f$ is $G$-equivariant. Furthermore, for
$r_1,r_2 \in R$,
\begin{equation*}
\begin{split}
\phi(g.(r_1x'r_2)) &= \phi(g.r_1)\phi(g.x')\phi(g.r_2) =
\phi(g.r_1)t\phi(g.x)\phi(g.r_2) \\ &=
t\phi(g.r_1)\phi(g.x)\phi(g.r_2) = t \phi(g.(r_1xr_2))\ .
\end{split}
\end{equation*}
This shows that $f$ is $(R,R)$-bilinear. Hence, defining $q$ to be
the class of $f$, we obtain the desired element $q\in \cC(R)^G$
mapping to our given $t \in \cC(S)^{\Htil}$, thereby proving
\eqref{E:orbC}.

It remains to construct $I$ as in the claim. Put
$$
D = \biggl( \bigcap_{\substack{x,y \in G\\ x^{-1}y \notin \Htil}}
x.(P:H) + y.(P:H) \biggr)^{[G:\Htil]-1}.
$$
Then $D$ is a $G$-stable ideal of $R$ satisfying $0 \neq \phi(D)$.
For the latter note that the finitely many ideals $x.(P:H) +
y.(P:H)$ are $H$-stable, since $H$ is normal, and none of them is
contained in $(P:H)$. By the Chinese remainder theorem
\cite[1.3]{kBmL96}, the image of the map
\begin{align*}
\mu \colon &R & &\into & &\prod_{g \in G/\Htil} R/g.(P:H) & &\iso & &\prod_{g \in G/\Htil} S \\
&\upin & & & &\qquad\upin & & & &\qquad\upin \\
 &\ r & &\mapsto & &(r + g.(P:H))_{g \in G/\Htil} & &\mapsto &
&(\phi(g^{-1}.r))_{g \in G/\Htil}
\end{align*}
contains the ideal $\prod_{g \in G/\Htil} \phi(D)$. Now put $I =
(\phi^{-1}(D_t):G)D$. This is certainly a $G$-stable ideal of $R$
satisfying $\phi(I) \subseteq D_t$\,. Suppose that $\phi(I) = 0$.
Since $\phi(D)$ is a nonzero $H$-stable ideal of the $H$-prime ring
$S$,  we must have
$$
(\phi^{-1}(D_t):G) = \bigcap_{g \in G/\Htil} g.\phi^{-1}(D_t)
\subseteq (P:H)
$$
and so $g.\phi^{-1}(D_t) \subseteq
(P:H)$ for some $g \in G$. But then $g.\phi^{-1}(D_t) \subsetneqq
\phi^{-1}(D_t)$ which is impossible because $\phi^{-1}(D_t)$ is
$\Htil$-stable and $\Htil$ has finite index in $G$. Therefore,
$\phi(I) \neq 0$. Finally, if $x \in I$ then $\phi(g.x) \in
D_t\phi(D)$ for all $g \in G$, and hence $t\phi(g.x) \in \phi(D)$.
Therefore, $(t\phi(g^{-1}.x))_{g \in G/\Htil} = \mu(x')$ for some
$x' \in R$, that is, $\phi(g^{-1}.x') = (t\phi(g^{-1}.x))$ holds for
all $g \in G/\Htil$. Since $\phi$ and $t$ are $\Htil$-invariant, it
follows that $\phi((gh)^{-1}.x') = (t\phi((gh)^{-1}.x))$ holds for
all $g \in G/\Htil$, $h \in \Htil$. Therefore, $\phi(g.x') =
t\phi(g.x)$ for all $g \in G$, as desired.

(b) This is just \eqref{E:orbC} with $H=1$.
\end{proof}



\subsection{$G$-rational ideals} \label{SS:GRat}

Assume now that $R$ is an algebra over some field $\k$, as in
Section~\ref{SS:rat}, and that $G$ acts on $R$ by $\k$-algebra
automorphisms.  A $G$-prime ideal $I$ of $R$ will be called
\emph{$G$-rational} if $\cC(R/I)^G = \k$. One can check as in
\S\ref{SSS:leftright} that the notion of $G$-rationality is
left-right symmetric.

Lemma~\ref{L:orbC}(a) with $H=1$ immediately implies the following

\begin{prop} \label{P:GRat}
The map $\Spec R \to \GSpec R$, $P \mapsto (P:G)$, in
Proposition~\ref{P:spec} restricts to a map $\Rat R \to \GRat R$.
\end{prop}

Unfortunately, the map $\Rat R \to \GRat R$ above need not be
surjective, even when the $G$-action on $R$ is locally ideal finite
in the sense of Proposition~\ref{P:spec}(b):

\begin{example}
Let $F \supset \k$ be any non-algebraic field extension satisfying
$F^G = \k$ for some subgroup $G$ of $\Gal(F/\k)$. For example, $F$
could be chosen to be the rational function field $\k(t)$ over an
infinite field $\k$ and $G = \k^*$ acting via $\lambda.f(t) =
f(\lambda^{-1}t)$ for $\lambda \in \k^*$. The $G$-action on $F$ is
clearly locally ideal finite and $\Qr(F) = \cC(F) = F$. Therefore,
the zero ideal of $F$ is $G$-rational, but $F$ has no rational
ideals.
\end{example}


\subsection{Algebras over a large algebraically closed base field}
\label{SS:affine}

We continue to assume that $R$ is an algebra over some field $\k$
and that $G$ acts on $R$ by $\k$-algebra automorphisms. The
following lemma is a version of \cite[Lemme 3.3]{cMrR81}.

\begin{lem} \label{L:affine}
Let $I \in \Spec R$ be given. Put $C = \cC(R/I)$ and consider the
natural map of $C$-algebras
$$
\psi \colon R_C = R \otimes_\k C \onto (R/I) \otimes_\k C \onto
(R/I)C
$$
where $(R/I)C \subseteq \Qr(R/I)$ is the central closure of $R/I$.
Then:
\begin{enumerate}
\item
$\Itil = \Ker \psi$ is a $C$-rational ideal of $R_C$.
\item
If $I \in \GRat R$ then, letting $G$ act on $R_C$ by $C$-linear
extension of the action on $R$, we have $$(\Itil:G)= I \otimes_\k C
\ .$$
\end{enumerate}
\end{lem}

\begin{proof}
Part (a) is clear, since $R_C/\Itil \cong (R/I)C$ is $C$-rational;
see \S\ref{SSS:rat1}.

For (b), note that the map $\psi$ is $G$-equivariant for the
diagonal $G$-action on $R_C = R \otimes_\k C$ and the usual
$G$-action on $(R/I)C \subseteq \Qr(R/I)$. Therefore, $\Itil$ is
stable under all automorphisms $g \otimes g$ with $g \in G$, and
hence we have $(g \otimes 1)(\Itil) = (1 \otimes g^{-1})(\Itil)$. We
conclude that
$$
(\Itil:G) = \bigcap_{g \in G} (1 \otimes g)(\Itil) = I \otimes_\k C\
,
$$
where the last equality uses the fact that $\Itil \cap R = I$ and
our hypothesis $C^G = \k$; see \cite[Cor. to Prop.~V.10.6]{nB81}.
\end{proof}

As an application of the lemma, we offer the following ``quick and
dirty" existence result for generic rational ideals.

\begin{prop} \label{P:affine}
Let $R$ be a countably generated algebra over an algebraically
closed base field $\k$ of infinite transcendence degree over its
prime subfield
and assume that the group $G$ is countably generated. Then every
prime ideal $I \in \GRat R$ has the form $I = (P:G)$ for some $P \in
\Rat R$.
\end{prop}

\begin{proof}
Let a prime $I \in \GRat R$ be given and let $\k_0$ denote the prime
subfield of $\k$.  By hypothesis on $R$, we have $\dim_{\k}R \le
\aleph_0$. Choosing a $\k$-basis $\cB$ of $R$ which contains a
$\k$-basis for $I$ and adjoining the structure constants of $R$ with
respect to $\cB$ to $\k_0$, we obtain a countable field $K$ with
$\k_0 \subseteq K \subseteq \k$. Putting $R_0 = \sum_{b \in \cB} Kb$
we obtain a $K$-subalgebra of $R$ such that $R = R_0\otimes_K \k$
and $I = I_0\otimes_K \k$, where $I_0 = I \cap R_0$. At the cost of
adjoining at most countably many further elements to $K$, we can
also make sure that $R_0$ is stable under the action of $G$. Thus,
$R_0/I_0$ is a $G$-stable $K$-subalgebra of $R/I$ and $R/I =
(R_0/I_0) \otimes_K \k$. Put $C = \cC(R_0/I_0)$ and note that
\eqref{E:scalar} implies that $C^G = K$, because $\cC(R/I)^G = \k$.
Thus, $I_0 \in \GRat R_0$ and Lemma~\ref{L:affine} yields an ideal
$\widetilde{I_0} \in \Rat(R_0 \otimes_K C)$ such that
$(\widetilde{I_0}:G) = I_0\otimes_K C$. Furthermore, since $R_0/I_0$
is countable, the field $C$ is countable as well; this follows from
Proposition~\ref{P:Qr}. By hypothesis on $\k$, there is a
$\k_0$-embedding of $C$ into $\k$; see \cite[Cor.~1 to
Th{\'e}or{\`e}me~V.14.5]{nB81}. Finally, Lemma~\ref{L:scalar}
implies that $P = \widetilde{I_0} \otimes_C \k$ is a rational ideal
of $(R_0 \otimes_K C) \otimes_C \k = R$ satisfying $(P:G) =
(I_0\otimes_K C) \otimes_C \k = I$, as desired.
\end{proof}


\section{Rational actions of algebraic groups} \label{S:algebraic}

In this section, we work over an algebraically closed base field
$\k$. Throughout, the group $G$ will be an affine algebraic group
over $\k$ and $R$ will be a $\k$-algebra on which $G$ acts by
$\k$-algebra automorphisms. The Hopf algebra of regular functions on
$G$ will be denoted by $\k[G]$. The notations introduced in
Section~\ref{S:Group} remain in effect. In addition, $\otimes$ will
stand for $\otimes_\k$.


\subsection{$G$-modules}  \label{SS:Gmod}

A $\k$-vector space $M$ is called a \emph{$G$-module} if there is a
linear representation
\begin{equation*} \label{E:Gmod1}
\rho_M \colon G \tto \GL(M)
\end{equation*}
satisfying
\begin{enumerate}
\item local finiteness: all $G$-orbits in $M$ generate
finite-dimensional subspaces of $M$, and
\item for every
finite-dimensional $G$-stable subspace $V \subseteq M$, the induced
group homomorphism $G \to \GL(V)$ is a homomorphism of algebraic
groups.
\end{enumerate}
As is well-known, these requirements are
equivalent to the existence of a $\k$-linear map
\begin{equation} \label{E:Gmod}
\Delta_M \colon M \tto M \otimes \k[G]
\end{equation}
which makes $M$ into an $\k[G]$-comodule; see Jantzen
\cite[2.7-2.8]{jcJ03} or Waterhouse \cite[3.1-3.2]{wW79} for
details. We will use the Sweedler notation
$$
\Delta_M(m) = \sum m_0 \otimes m_1 \qquad (m \in M)
$$
as in Montgomery \cite{sM93}. Writing $\rho_M(g)(m) = g.m $, we have
\begin{equation} \label{E:action}
g.m  = \sum m_0 m_1(g) \qquad (g \in G, m \in M)\ .
\end{equation}

Linear representations $\rho_M$ as above are often called
\emph{rational}. Tensor products of rational representations of $G$
are again rational, and similarly for sums, subrepresentations and
homomorphic images of rational representations.

\begin{example}
If the group $G$ is finite then $G$-modules are the same as (left)
modules $M$ over the group algebra $\k G$ and all linear
representations of $G$ are rational. Indeed, in this case, $\k[G]$
is the linear dual of $\k G$, that is, the $\k$-vector space of all
functions $G \to \k$ with pointwise addition and multiplication. The
map $\Delta_M \colon M \to M \otimes \k[G]$ is given by
$$
\Delta_M(m) = \sum_{x \in G} x.m \otimes p_x\ ,
$$
where $p_x \in \k[G] = (\k G)^*$ is defined by $p_x(y) =
\delta_{x,y}$ (Kronecker delta) for $x,y \in G$.
\end{example}


\subsection{Some properties of $G$-modules}  \label{SS:Gmod2}

Let $M$ be a $G$-module. The coaction $\Delta_M$ in \eqref{E:Gmod}
%
%
is injective. In fact, extending $\Delta_M$ to a map
\begin{equation} \label{E:Gmod2}
\Delta_M \colon M \otimes \k[G] \tto M \otimes \k[G]
\end{equation}
by $\k[G]$-linearity, we obtain an \emph{automorphism} of $M \otimes
\k[G]$: the inverse of $\Delta_M$ is the $\k[G]$-linear extension of
the map
$(\Id_M \otimes S) \circ \Delta_M \colon M  \tto M \otimes \k[G]$,
where $S \colon \k[G] \to \k[G]$  is the antipode of $\k[G]$:
$(Sf)(g) = f(g^{-1})$ for $g \in G$.

Furthermore, $G$-stable cores can be computed with $\Delta_M$ as
follows.

\begin{lem} \label{L:Gcore}
For any $\k$-subspace $V$ of $M$, we have $(V:G) = \Delta_M^{-1}(V
\otimes \k[G])$
\end{lem}

\begin{proof}
Fix a $\k$-basis $\{ v_i \}$ of $V$ and let $\{ w_j \}$ be a
$\k$-basis for a complement of $V$ in $M$. For $m \in M$, we have
$\Delta_M(m) = \sum_i v_i \otimes f_i + \sum_j w_j \otimes h_j$ with
uniquely determined $f_i, h_j \in \k[G]$. Moreover,
$$
\Delta_M(m) \in V \otimes \k[G] \iff \text{all $h_j = 0$} \iff
\forall g \in G : g.m = \sum_i v_i f_i(g) \in V\ .
$$
This proves the lemma.
\end{proof}


\subsection{Regular representations and intertwining formulas}
\label{SS:regular}

The right and left \emph{regular representations} of $G$ are defined
by
\begin{equation} \label{E:reg}
\begin{aligned}
&\rho_r \colon G \tto \GL(\k[G])\,,\quad &
&\left(\rho_r(x)f\right)(y)
= f(yx)\ , \\
&\rho_\ell \colon G \tto \GL(\k[G])\,,\quad &
&\left(\rho_\ell(x)f\right)(y) = f(x^{-1}y)\ .
\end{aligned}
\end{equation}
for $x,y \in G$.
Both regular representations are rational. The right regular
representation comes from the comultiplication $\Delta \colon \k[G]
\to \k[G] \otimes \k[G]$ of the Hopf algebra $\k[G]$; in the usual
Sweedler notation, it is given by $\Delta f = \sum f_1 \otimes f_2$,
where $f(xy) = \sum f_1(x)f_2(y)$ for $x,y \in G$. Similarly, the
left regular representation comes from $\left( S \otimes \Id_{\k[G]}
\right) \circ \Delta \circ S \colon \k[G] \to \k[G] \otimes \k[G]$.

Now let $M$ be a $G$-module. Then the rational representations $1_M
\otimes \rho_\ell \colon G \to \GL(M \otimes \k[G])$ and $\rho_M
\otimes \rho_\ell \colon G \to \GL(M \otimes \k[G])$ are intertwined
by the automorphism $\Delta_M$ of \eqref{E:Gmod2}: for all $g \in
G$, we have
\begin{equation} \label{E:reg1}
\Delta_M \circ \left( 1_M \otimes \rho_\ell \right)(g) = \left(
\rho_M \otimes \rho_\ell \right)(g) \circ \Delta_M \ .
\end{equation}
Similarly,
\begin{equation} \label{E:reg2}
\Delta_M \circ \left( \rho_M \otimes \rho_r \right)(g)  = \left( 1_M
\otimes \rho_r \right)(g) \circ \Delta_M \ .
\end{equation}
To prove \eqref{E:reg2}, for example, one checks that both sides of
the equation send $m \otimes f \in M \otimes \k[G]$ to the function
$G \to M$, $x \mapsto xg.m f(xg)$.


\subsection{Rational group actions}  \label{SS:rational}

The action of $G$ on the $\k$-algebra $R$ is said to be
\emph{rational} if it makes $R$ a $G$-module in the above sense. The
map
$$
\Delta_R \colon R \to  R \otimes \k[G]
$$
is then a map of $\k$-algebras; equivalently, $R$ is a right
$\k[G]$-comodule algebra. Since rational actions are locally finite,
they are certainly locally ideal finite in the sense of
Proposition~\ref{P:spec}(b). Therefore, the $G$-primes of $R$ are
exactly the ideals of $R$ of the form $(P:G)$ for $P \in \Spec R$.
In particular, $G$-prime ideals of $R$ are semiprime; for a more
precise statement, see Corollary~\ref{C:connected} below. Moreover,
the $\k[G]$-linear extension of $\Delta_R$ is an automorphism of
$\k[G]$-algebras
\begin{equation} \label{E:Delta}
\Delta_R \colon R \otimes \k[G] \iso  R \otimes \k[G] \ .
\end{equation}

We now consider the extended $G$-action on the Amitsur-Martindale
ring of quotients $\Qr(R)$; see \S\ref{SS:GSpecC}. This action is
usually not rational, even if $G$ acts rationally on $R$. Part (b)
of the following lemma, for classical quotient rings of semiprime
Goldie rings, is due to M{\oe}glin and Rentschler
\cite[I.22]{cMrR86}.

\begin{lem} \label{L:closed}
Assume that $G$ acts rationally on $R$. Then:
\begin{enumerate}
\item
The centralizer $C_G(T) = \{ g \in G \mid g.q = q \ \forall q \in T
\}$ of every subset $T \subseteq \Qr(R)$ is a closed subgroup of
$G$.
\item
Let $V \subseteq \Qr(R)$ be a $G$-stable $\k$-subspace of $\Qr(R)$.
The $G$-action on $V$ is rational if and only if it is locally
finite.
\end{enumerate}
\end{lem}

\begin{proof}
(a) In view of Proposition~\ref{P:Qr}(iii), the condition for an
element $g \in G$ to belong to $C_G(T)$ can be stated as
$$
\forall q \in T, r \in D_q : (q - g.q)g.r = 0 \ ,
$$
where $D_q$ is as in \eqref{E:Dq}. Using the notation of
\eqref{E:action}, we have
$$
(q - g.q)g.r = q(g.r) - g.(qr) = \sum qr_0 r_1(g) - \sum
(qr)_0(qr)_1(g) \ .
$$
Thus, putting $f_{r,q} = \sum qr_0 \otimes r_1 - \sum (qr)_0 \otimes
(qr)_1 \in \Qr(R) \otimes \k[G]$, we see that $g \in C_G(T)$ if and
only if $f_{r,q}(g) = 0$ holds for all $q \in T$ and all $r \in
D_q$. Since each equation $f_{r,q}(g) = 0$ defines a closed subset
of $G$, part (a) follows.

(b) Necessity is clear. So assume that the $G$-action on $V$ is
locally finite. Put $S = R \otimes \k[G]$ and consider the
$\k[G]$-algebra automorphism $\Delta_R \in \Aut(S)$ as in
\eqref{E:Delta} and its extension $\Delta \in \Aut(\Qr(S))$. We must
show that, under the canonical embedding $\Qr(R) \into \Qr(S)$ as in
\S\ref{SS:scalar}, we have
\begin{equation} \label{E:closed1}
\Delta(V) \subseteq V \otimes \k[G] \ .
\end{equation}
Since the action of $G$ on $V$ is locally finite, we may assume that
$V$ is finite-dimensional. Therefore, the ideal $D_V = \bigcap_{q
\in V} D_q$ belongs to $\sE(R)$ and $D_V$ is $G$-stable, since $V$
is. Lemma~\ref{L:Gcore} implies that $\Delta(D_V \otimes \k[G]) =
D_V \otimes \k[G]$, and hence
$$
\Delta(V) (D_V \otimes \k[G]) = \Delta(V (D_V \otimes \k[G]))
\subseteq S \ .
$$
This shows that the subspace $\Delta(V) \subseteq \Qr(S)$ actually
is contained in $\Qr(R) \otimes \k[G]$, and \eqref{E:closed1}
follows from Lemma~\ref{L:Gcore}, since $V = (V:G)$.
\end{proof}

\emph{From now on, the $G$-action on $R$ is understood to be
rational.}


\subsection{Connected groups}
\label{SS:connected}

The group $G$ is \emph{connected} if and only if the algebra $\k[G]$
is a domain. In this case,
$$
\k(G) = \Fract \k[G]
$$
will denote the field of rational functions on $G$. The group $G$
acts on $\k(G)$ by the natural extensions of the right and left
regular actions $\rho_r$ and $\rho_\ell$ on $\k[G]$; see
\S\ref{SS:regular}.

Part (a) of the following result is due to Chin \cite[Corollary
1.3]{wC92}; the proof given below has been extracted from
\cite[3.6]{nVE98}. The proof of part (c) follows the outline of the
arguments in \cite[I.29, $2^e$ {\'e}tape]{cMrR86}.

\begin{prop} \label{P:connected}
Assume that $G$ is connected. Then:
\begin{enumerate}
\item
$(P:G)$ is prime for every $P \in \Spec R$. Therefore, the
$G$-primes of $R$ are exactly the $G$-stable primes of $R$.
\item
Assume that $R$ is prime and every nonzero ideal $I$ of $R$
satisfies $(I:G) \neq 0$. Then $G$ acts trivially on $\cC(R)$.
\item
If $R$ is $G$-rational then the field extension $\cC(R)/\k$ is
finitely generated. In fact, there is a $G$-equivariant
$\k$-embedding of fields $\cC(R) \into \k(G)$, with $G$ acting on
$\k(G)$ via the right regular representation $\rho_r$.
\end{enumerate}
\end{prop}

\begin{proof}
(a) It suffices to show that $(P:G)$ is prime for each prime $P$;
the last assertion is then a consequence of
Proposition~\ref{P:spec}.

By \S\ref{SS:rational}, we know that the homomorphism $\Delta_R
\colon R \to R \otimes \k[G]$ is centralizing. Therefore, there is a
map
$\Spec( R \otimes \k[G]) \to \Spec R$, $Q \mapsto \Delta_R^{-1}(Q)$.
In view of Lemma~\ref{L:Gcore}, it therefore suffices to show that
$P \otimes \k[G]$ is prime whenever $P$ is. But the algebra $\k[G]$
is contained in some finitely generated purely transcendental field
extension $F$ of $\k$; see Borel \cite[18.2]{aB91}.  Thus, we have a
centralizing extension of algebras
$$
(R/P) \otimes \k[G] \subseteq   (R/P) \otimes F \ .
$$
Since $(R/P) \otimes F$ is clearly prime, $(R/P) \otimes \k[G]$ is
prime as well as desired.

(b) We first prove the following special case of (b) which is
well-known; see \cite[Prop.~A.1]{nVE93}.
\begin{claim}
If $R$ is a field then $G$ acts trivially on $R$.
\end{claim}
Since $G$ is the union of its Borel subgroups (\cite[11.10]{aB91}),
we may assume that $G$ is solvable. Arguing by induction on a
composition series of $G$ (\cite[15.1]{aB91}), we may further assume
that $G$ is the additive group $\Ga$ or the multiplicative group
$\Gm$. Therefore, $R \otimes \k[G]$ is a polynomial algebra or a
Laurent polynomial algebra over $R$. In either case, $R$ is the
unique largest subfield of $R \otimes \k[G]$, because $R \otimes
\k[G]$ has only ``trivial" units: the nonzero elements of $R$ if $R
\otimes \k[G] = R[t]$, and the elements of the form $r t^m$ with $0
\neq r \in R$ and $m \in \ZZ$ if $R \otimes \k[G] = R[t^{\pm 1}]$.
Consequently, the map $\Delta_R \colon R \to R \otimes \k[G]$ has
image in $R \otimes 1$ which in turn says that $G$ acts trivially on
$R$. This proves the Claim.

Now let $R$ be a prime $\k$-algebra such that $(I:G)$ is nonzero for
every nonzero ideal $I$ of $R$. By the Claim, it suffices to show
that the $G$-action on $\cC(R)$ is rational, and by
Lemma~\ref{L:closed} this amounts to showing that $G$-action on
$\cC(R)$ is locally finite. So let $q \in \cC(R)$ be given and
consider the ideal $D_q$ of $R$ as in \eqref{E:Dq}. By hypothesis,
we may pick a nonzero element $d \in (D_q:G)$. The $G$-orbit $G.d$
generates a finite-dimensional $\k$-subspace $V \subseteq D_q$.
Moreover, $qV$ is contained in a finite-dimensional $G$-stable
subspace $W \subseteq R$. Therefore, for all $g,h \in G$, we have
$(g.q)(h.d) = g.(q(g^{-1}h.d)) \in W$, and hence $QV \subseteq W$,
where $Q \subseteq \cC(R)$ denotes the $\k$-subspace that is
generated by the orbit $G.q$. Thus, multiplication gives a linear
map $Q \to \Hom_{\k}(V,W)$ which is injective, because $V \neq 0$
and nonzero elements of $\cC(R)$ have zero annihilator in $R$. This
shows that $Q$ is finite-dimensional as desired.

(c) Put $C = \cC(R)$ and $K = \k(G)$, the field of rational
functions on $G$, that is, the field of fractions of the algebra
$\k[G]$. The algebra $R_K = R \otimes K$ is prime by (a) and its
proof, and by \eqref{E:scalar} there is a tower of fields
\begin{equation*}
C \into \Fract(C \otimes K) \into \cC(R_K) \ .
\end{equation*}
We will first show that $C$ is a finitely generated field extension
of $\k$. Since $K/\k$ is finitely generated, the field $\Fract(C
\otimes K)$ is certainly finitely generated over $C$. Thus, it will
suffice to construct a $C$-algebra embedding $C \otimes C \into
\Fract(C \otimes K)$.

To construct such an embedding, consider the natural epimorphism of
$\cC(R_K)$-algebras $RC \otimes_C \cC(R_K) \onto R_K\cC(R_K)$. By
Lemma~\ref{L:Central}(b), this map is injective, because it is
clearly injective on $RC$. Thus,
\begin{equation} \label{E:iso}
RC \otimes_C \cC(R_K) \iso R_K\cC(R_K) \ .
\end{equation}
Let $\delta$ be the $K$-algebra automorphism of $R_K$ that is
defined by $K$-linear extension of the $G$-coaction $\Delta_R \colon
R \otimes \k[G] \iso R \otimes \k[G]$ in \eqref{E:Delta}:
\begin{equation} \label{E:delta}
\delta = \Delta_R \otimes_{\k[G]}\Id_K \colon R_K \iso R_K \ .
\end{equation}
Let $\til{\delta}$ be the unique extension of $\delta$ to an
automorphism of the central closure $R_K\cC(R_K)$ of $R_K$. Clearly,
$\til{\delta}$ sends the $\cC(R_K) = \cen(R_K\cC(R_K))$ to itself.
We claim that
\begin{equation} \label{E:deltaC}
\til{\delta}(C) \subseteq \Fract(C \otimes K)\ ;
\end{equation}
so $\til{\delta}$ also sends the $\Fract(C \otimes K)$ to itself. In
order to see this, pick $q \in C$ and $d \in D_q$. Then
$$
\til{\delta}(q)\Delta_R(d) = \til{\delta}(q)\til{\delta}(d) =
\til{\delta}(qd) = \Delta_R(qd)
$$
holds in $R_K\cC(R_K)$. Here, both $\Delta_R(qd)$ and $\Delta_R(d)$
belong to $R_K \subseteq RC \otimes_C (C \otimes K)$. Fixing a
$C$-basis $B$ for $RC$ and writing $\Delta_R(qd) =  \sum_{b \in B} b
x_b$ and $\Delta_R(d)  =  \sum_{b \in B} b y_b$ with $x_b,y_b \in C
\otimes K$, the above equation becomes
$$
\sum_{b \in B} b \til{\delta}(q)y_b = \sum_{b \in B} b x_b \ .
$$
Now \eqref{E:iso} yields $\til{\delta}(q)y_b = x_b$ for all $b$,
which proves \eqref{E:deltaC}. Now, for the desired embedding,
consider the $C$-algebra map
\begin{equation} \label{E:mu}
\mu \colon C \otimes C \tto \Fract(C \otimes K)\ , \quad c \otimes
c' \mapsto c \til{\delta}(c') \ .
\end{equation}
We wish to show that $\mu$ is injective. To this end, note that the
$G$-action $\rho_R$ on $R$ extends uniquely to an action $\rho_{RC}$
on the central closure $RC$, and the $G$-action $1_R \otimes \rho_r$
on $R_K$ extends uniquely to the central closure $R_K\cC(R_K)$.
Denoting this latter action by $\til{\rho_r}$, the intertwining
formula \eqref{E:reg2} implies that $\til{\delta} \circ \rho_{RC}(g)
= \til{\rho_r}(g) \circ \til{\delta} \colon RC \to R_K\cC(R_K)$ for
all $g \in G$. This yields
\begin{equation} \label{E:mu1}
\mu \circ (\Id_C \otimes \rho_C(g)) = \til{\rho_r}(g) \circ \mu
\end{equation}
for all $g \in G$. Thus, the ideal $\Ker\mu$ of $C \otimes C$ is
stable under $(1_C \otimes \rho_C)(G)$. Finally, since $C^G = \k$,
we may invoke \cite[Cor. to Prop.~V.10.6]{nB81} to conclude that
$\Ker\mu$ is generated by its intersection with $C \otimes 1$, which
is zero. This shows that $\mu$ is injective, and hence the field
extension $C/\k$ is finitely generated.

It remains to construct a $G$-equivariant embedding $C \into K$,
with $G$ acting on $\k(G)$ via the right regular representation
$\rho_r$ as above. For this, we specialize \eqref{E:mu} as follows.
Write $C = \Fract A$ for some affine $\k$-subalgebra $A \subseteq
C$. Then $\Fract(C \otimes K) = \Fract(A \otimes \k[G])$, and hence
$$
\mu(A \otimes A) \subseteq (A \otimes \k[G])[s^{-1}]
$$
for some $0 \neq s \in A \otimes \k[G]$. By generic flatness (e.g.,
Dixmier \cite[2.6.3]{jD96}), there further exists $0 \neq f \in A
\otimes A$ so that $(A \otimes \k[G])[\mu(f)^{-1}s^{-1}]$ is free
over $(A \otimes A)[f^{-1}]$ via $\mu$. Now choose some maximal
ideal $\fm$ of $A$ with $f \notin \fm \otimes A$. Let $\bar{f}$
denote the image of $f$ in $(A \otimes A)/(\fm \otimes A) \cong A$,
and let $\bar{s}$ denote the image of $s$ in $(A \otimes \k[G])/(\fm
\otimes \k[G]) \cong \k[G]$. Since $\mu(\fm \otimes A) = \fm\mu(A
\otimes A)$, the map $\mu\vert_{A\otimes A}$ passes down to a map
$$
\bar{\mu} \colon A[\bar{f}^{-1}] \tto B:=
\k[G][\bar{\mu}(\bar{f})^{-1}\bar{s}^{-1}]
$$
making $B$ a free $A[\bar{f}^{-1}]$-module. Consequently,
$\bar{\mu}$ extends uniquely to an embedding of the fields of
fractions, $\Fract A[\bar{f}^{-1}] = C \into \Fract B = K$. Finally,
\eqref{E:mu1} implies that this embedding is $G$-equivariant, which
completes the proof of (c).
\end{proof}

Returning to the case of a general affine algebraic group $G$, we
have the following

\begin{cor} \label{C:connected}
Every $I \in \GSpec R$ has the form $I = (Q:G)$ for some $Q \in
\Spec R$ with $[G:G_Q] < \infty$. Moreover, $\cC(I)^G \cong
\cC(Q)^{G_Q}$.
\end{cor}

\begin{proof}
We know that $I = (P:G)$ for some $P \in \Spec R$; see
\S\ref{SS:rational}. Let $G^0$ denote the connected component of the
identity in $G$; this is a connected normal subgroup of finite index
in $G$ (e.g., Borel \cite[1.2]{aB91}). Put $Q = (P:G^0)$. Then
Proposition~\ref{P:connected}(a) tells us that $Q$ is prime.
Furthermore, $I = (Q:G)$ and $G^0 \subseteq G_Q$; so $[G:G_Q] <
\infty$. The isomorphism $\cC(I)^G \cong \cC(Q)^{G_Q}$ follows from
Lemma~\ref{L:orbC}(b).
\end{proof}


\subsection{The fibres of the map \eqref{E:ratmap}} \label{SS:fibres}

Assume that $G$ is connected. Our next goal is to give a description
of the fibres of the map $\Rat R \to \GRat R$, $P \mapsto (P:G)$ in
Proposition~\ref{P:GRat}. Following \cite{kBkG02} we denote the
fibre over a given $I \in \GRat R$ by $\Rat_IR$:
$$
\Rat_IR = \{ P \in \Rat R \mid (P:G) = I \} \ .
$$
The group $G$ acts on $\Rat_IR$ via the given action $\rho_R$ on
$R$.

Recall that the group $G$ acts on the rational function field
$\k(G)$ by the natural extensions of the regular representations
$\rho_r$ and $\rho_\ell$. We let
$$
\Hom_G(\cC(R/I),\k(G))
$$
denote the collection of all $G$-equivariant $\k$-algebra
homomorphisms $\cC(R/I) \to \k(G)$ with $G$ acting on $\k(G)$ via
the right regular action $\rho_r$. The left regular action
$\rho_\ell$ of $G$ on $\k(G)$ yields a $G$-action on the set
$\Hom_G(\cC(R/I),\k(G))$.

\begin{thm} \label{T:fibres}
Let $I \in \GRat R$ be given. There is a $G$-equivariant bijection
$$
\Rat_IR \tto \Hom_G(\cC(R/I),\k(G))\ .
$$
\end{thm}

\begin{proof}
Replacing $R$ by $R/I$, we may assume that $I=0$. In particular, $R$
is prime by Proposition~\ref{P:connected}. We will also put $C =
\cC(R)$ and $K = \k(G)$ for brevity. For every $P \in \Rat R$ with
$(P:G) = 0$, we will construct an embedding of fields
$$
\psi_P \colon C \into K
$$
such that the following hold:
\begin{enumerate}
\item
$\psi_P(g.c) = \rho_r(g)(\psi_P(c))$ and $\psi_{g.P} =
\rho_\ell(g)\circ \psi_P$ holds for all $g \in G$, $c \in C$;
\item
if $P,Q \in \Rat R$ are such that $(Q:G) = (P:G) = 0$ but $Q \neq P$
then $\psi_Q \neq \psi_P$;
\item
given a $G$-equivariant embedding $\psi \colon C \into K$, with $G$
acting on $K$ via $\rho_r$, we have $\psi = \psi_P$ for some $P \in
\Rat R$ with $(P:G) = 0$.
\end{enumerate}
This will prove the theorem.

In order to construct $\psi_P$, consider the $K$-algebra $(R/P)_K =
(R/P)\otimes K$. This algebra is rational by Lemma~\ref{L:scalar}.
We have a centralizing $\k$-algebra homomorphism
\begin{equation} \label{E:phiP}
\phi_P \colon R \stackrel{\Delta_R}{\tto} R \otimes \k[G]
 \stackrel{\text{can.}}{\tto}
(R/P)_K \ ,
\end{equation}
where the canonical map $R \otimes \k[G] \to (R/P)_K$ comes from the
embedding $\k[G] \into K$ and the epimorphism $R \onto R/P$. Since
$(P:G)=0$, Lemma~\ref{L:Gcore} implies that $\phi_P$ is injective.
Since $(R/P)_K$ is prime, it follows that $\cC_{\phi_P} = C$ holds
in Lemma~\ref{L:C}. Hence $\phi_P$ extends uniquely to a
centralizing $\k$-algebra monomorphism
\begin{equation} \label{E:phiP2}
\til{\phi}_P \colon RC \into (R/P)_K\cC((R/P)_K) = (R/P)_K
\end{equation}
sending $C$ to $\cC((R/P)_K) = K$. Thus we may define $\psi_P:=
\til{\phi}_P\vert_{C} \colon C \into K$. It remains to verify
properties (a) - (c).

Part (a) is a consequence of the intertwining formulas
\eqref{E:reg1} and \eqref{E:reg2}. Indeed,  \eqref{E:reg2} implies
that
$\phi_P(g.r) = \rho_r(g)(\phi_P(r))$
holds for all $g \in G$ and $r \in R$. In view of
Proposition~\ref{P:Qr}(ii), this identity is in fact valid for
$\til{\phi}_P$ and all $r \in RC$, which proves the first of the
asserted formulas for $\psi_P$ in (a). For the second formula,
consider the map $(\phi_P)_K$ that is defined by $K$-linear
extension of \eqref{E:phiP} to $R_K = R \otimes K$; this is the
composite
\begin{equation} \label{E:phiP3}
(\phi_P)_K  \colon R_K \stackrel{\delta}{\tto} R_K
\stackrel{\text{can.}}{\onto} (R/P)_K \ ,
\end{equation}
where $\delta$ is as in \eqref{E:delta}. The map $(\rho_R \otimes
\rho_\ell)(g)$ gives ring isomorphisms $R_K \iso R_K$ and $(R/P)_K
\iso (R/g.P)_K$ such that the following diagram commutes:
$$
\xymatrix{%
R_K \ar[r]^-{\sim} \ar@{>>}[d]_{\text{\rm can.}} &
R_K\ar@{>>}[d]^{\text{\rm can.}}
\\
(R/P)_K \ar[r]^-{\sim} & (R/g.P)_K }
$$
The intertwining formula \eqref{E:reg1} implies that, for all $g \in
G$,
$$
(\phi_{g.P})_K \circ (1_R \otimes \rho_\ell)(g) = (\rho_R \otimes
\rho_\ell)(g) \circ (\phi_P)_K\ .
$$
Restricting to $R$ we obtain $\phi_{g.P} = (\rho_R \otimes
\rho_\ell)(g) \circ \phi_P$, and this becomes $\psi_{g.P} =
\rho_\ell(g)\circ \psi_P$ on $C$. This finishes the proof of (a).

For (b), let
\begin{equation} \label{E:phiP4}
(\til{\phi}_P)_K  \colon (RC)_K = RC \otimes K  \onto (R/P)_K
\end{equation}
be defined by $K$-linear extension of \eqref{E:phiP2} and put
$\til{P} = \Ker(\til{\phi}_P)_K$. Let $Q \in \Rat R$ be given such
that $(Q:G) = 0$ and let $\til{Q} = \Ker(\til{\phi}_Q)_K$ be defined
analogously. If $Q \neq P$ then $\til{Q}$ and $\til{P}$ and are
distinct primes of $(RC)_K$; in fact, $\til{Q} \cap R_K \neq \til{P}
\cap R_K$, because the restriction of $(\til{\phi}_P)_K $ to $R_K$
is given by \eqref{E:phiP3}. Since both $\til{Q}$ and $\til{P}$ are
disjoint from $RC$, Lemma~\ref{L:Central}(c) gives $\til{P} \cap C_K
\neq \til{Q} \cap C_K$. This shows that $(\psi_P)_K$ and
$(\psi_Q)_K$ have distinct kernels, and so $\psi_P \neq \psi_Q$
proving (b).

Finally, for (c), let $\psi \colon C \into K$ be some
$G$-equivariant embedding. Define a $K$-algebra map
\begin{equation*}
\Psi \colon R_K \tto S = RC \otimes_C K
\end{equation*}
by $K$-linear extension of the canonical embedding $R \into RC$.
Note that, for $c \in C$,
\begin{equation} \label{E:Psi}
c \otimes 1 = 1 \otimes \psi(c)
\end{equation}
holds in $S$. Put
$$
P = \delta(\Ker\Psi) \cap R \ ,
$$
with $\delta$ as in \eqref{E:delta}. We will show that $P$ is the
desired rational ideal.

The algebra $S$ is $K$-rational, by Lemma~\ref{L:scalar}, and $G$
acts on $S$ via $\rho_{RC} \otimes_C \rho_r$, where $\rho_{RC}$ is
the unique extension of the $G$-action $\rho_R$ from $R$ to the
central closure $RC$. The map $\Psi$ is $G$-equivariant for this
action and the diagonal $G$-action $\rho_R \otimes \rho_r$ on $R_K$.
Furthermore, by \eqref{E:reg2}, the automorphism $\delta^{-1} \colon
R_K \iso R_K$ is equivariant with respect to the $G$-actions $1_R
\otimes \rho_r$ on the first copy of $R_K$ and $\rho_R \otimes
\rho_r$ on the second $R_K$. Therefore, the composite
$\Psi\circ\delta^{-1} \colon R_K \to S$ is equivariant for the
$G$-actions $1_R \otimes \rho_r$ on $R_K$ and $\rho_{RC} \otimes_C
\rho_r$ on $S$. Now consider the centralizing monomorphism of
$\k$-algebras
$$
\mu \colon R/P \into R_K/\delta(\Ker\Psi)
\underset{\delta^{-1}}{\iso} R_K/\Ker\Psi \underset{\Psi}{\into} S \
.
$$
By the foregoing, we have $\mu(R/P) \subseteq S^G$, the
$\k$-subalgebra of $G$-invariants in $S$. Since $S$ is prime, we
have $\cC_\mu = \cC(R/P)$ in Lemma~\ref{L:C}. Hence, $\mu$ extends
uniquely to a monomorphism $\til{\mu} \colon R/P\cC(R/P) \into
S\cC(S) = S$ sending $\cC(R/P)$ to $\cC(S) = K$. Therefore,
$\til{\mu}(\cC(R/P)) \subseteq K^G = \k$, which proves that $P$ is
rational. Furthermore, by Lemma~\ref{L:Gcore}, we have $(P:G) =
\Delta_R^{-1}(P \otimes \k[G]) \subseteq
\delta^{-1}(\delta(\Ker\Psi)) = \Ker\Psi$. Since $\Psi$ is mono on
$R$, we conclude that $(P:G) = 0$. It remains to show that $\psi =
\psi_P$. For this, consider the map $\til{\phi}_P$ of
\eqref{E:phiP2}; so $\psi_P = \til{\phi}_P\vert_{C}$. For $q \in C$,
$d \in D_q$ we have
$$
\delta(qd) \bmod P \otimes K = \til{\phi}_P(qd) =
\til{\phi}_P(q)\til{\phi}_P(d) = \delta(\psi_P(q)d) \bmod P \otimes
K
$$
because $\psi_P(q) \in K$ and $\delta$ is $K$-linear. It follows
that $\psi_P(q)d - qd \in \Ker\Psi$; so $0 = \psi_P(q)\Psi(d) -
\Psi(qd) = qd \otimes_C 1 = \psi(q)\Psi(d)$, where the last equality
holds by \eqref{E:Psi}. This shows that $\psi_P(q) = \psi(q)$,
thereby completing the proof of the theorem.
\end{proof}


\subsection{Proof of Theorem~\ref{T:goal}}
\label{SS:proof}

We have to prove:
\begin{enumerate}
\item[(1)] given $I \in \GRat R$, there is a $P
\in \Rat R$ such that $I = (P:G)$;
\item[(2)] if $P,P' \in \Rat R$ satisfy $(P:G) = (P':G)$ then $P' = g.P$
for some $g \in G$.
\end{enumerate}

\subsubsection{} \label{SSS:G0reduction}
We first show that it suffices to deal with the case of connected
groups. Let $G^0$ denote the connected component of the identity in
$G$, as before, and assume that both (1) and (2) hold for $G^0$.

In order to prove (1) for $G$, let $I \in \GRat R$ be given. By
Corollary~\ref{C:connected}, there exists $Q \in \Spec R$ with $I =
(Q:G)$,  $G^0 \subseteq G_Q$ and $\cC(R/Q)^{G_Q} = \k$. Since
$G_Q/G^0$ is finite, it follows that $Q$ is in fact $G^0$-rational.
Inasmuch as (1) holds for $G^0$, there exists $P \in \Rat R$ with $Q
= (P:G^0)$. It follows that $(P:G) = (Q:G) = I$, proving (1).

Now suppose that $(P:G) = (P':G)$ for $P,P' \in \Rat R$. Putting
$P^0 = (P:G^0)$ we have $(P:G) = \bigcap_{x \in G/G^0} x.P^0 =
\bigcap_{x \in G/G^0} (x.P:G^0)$, a finite intersection of
$G^0$-prime ideals of $R$. Similarly for $P'^0 = (P':G^0)$. The
equality $(P:G) = (P':G)$ implies that $(P':G^0) = (x.P:G^0)$ for
some $x \in G$. (Note that if $V \subseteq g.V$ holds for some
$\k$-subspace $V \subseteq R$ and some $g \in G$ then we must have
$V = g.V$, because the $G$-action on $R$ is locally finite.)
Invoking (2) for $G^0$, we see that $P' = yx.P$ for some $y \in
G^0$, which proves (2) for $G$.

\subsubsection{} \label{SS:proof2}
Now assume that $G$ is connected. In view of Theorem~\ref{T:fibres},
proving (1) amounts to showing that there is a $G$-equivariant
$\k$-algebra homomorphism $\cC(R/I) \to \k(G)$ with $G$ acting on
$\k(G)$ via the right regular action $\rho_r$. But this has been
done in Proposition~\ref{P:connected}(c). For part (2), it suffices
to invoke Theorem~\ref{T:fibres} in conjunction with the following
result which is the special case of Vonessen \cite[Theorem
4.7]{nVE98} for connected $G$.

\begin{prop} \label{P:injectivity}
Let $G$ act on $\k(G)$ via $\rho_r$ and let $F$ be a $G$-stable
subfield of $\k(G)$ containing $\k$. Let $\Hom_G(F,\k(G))$ denote
the collection of all $G$-equivariant $\k$-algebra homomorphisms
$\phi \colon F \to \k(G)$. Then the $G$-action on $\Hom_G(F,\k(G))$
that is given by $g.\phi = \rho_\ell(g)\circ \phi$ is transitive.
\end{prop}

This completes the proof of Theorem~\ref{T:goal}. \qed

\subsubsection{} \label{SSS:surjectivity}
It is tempting to try and prove (1) above in the following more
direct fashion. Assume that $R$ is $G$-prime and choose an ideal $P$
of $R$ that is maximal subject to the condition $(P:G) = 0$. This is
possible by the proof of Proposition~\ref{P:spec}(b) and we have
also seen that $P$ is prime. I don't know if the ideal $P$ is
actually rational. This would follow if the field extension
$\cC(R)^G \into \cC(R/P)^{G_P}$ in Lemma~\ref{L:orbC} were algebraic
in the present situation. Indeed, every ideal $I$ of $R$ with $I
\supsetneq P$ satisfies $(I:G) \neq 0$, and hence $(I:H) \supsetneq
P$. Therefore, Proposition~\ref{P:connected}(b) tells us that the
connected component of the identity of $G_P$ acts trivially on
$\cC(R/P)$ and so $\cC(R/P)$ is finite over $\cC(R/P)^{G_P}$.


\begin{ack}
It is a pleasure to thank Rudolf Rentschler for his comments,
questions, and guidance to the literature on enveloping algebras.
The author also thanks Ken Goodearl, Jens Carsten Jantzen, Peter
Linnell, and Lance Small for helpful comments at various stages of
this work, and the referee for a careful reading of the article and
valuable suggestions.
\end{ack}


\bibliographystyle{amsplain}

\def\cprime{$'$}
\providecommand{\bysame}{\leavevmode\hbox
to3em{\hrulefill}\thinspace}
\providecommand{\MR}{\relax\ifhmode\unskip\space\fi MR }
\providecommand{\MRhref}[2]{%
  \href{http://www.ams.org/mathscinet-getitem?mr=#1}{#2}
} \providecommand{\href}[2]{#2}


\end{document}